\documentclass[12pt]{article}

\usepackage{amsmath,amssymb,enumerate,bbm}
\newcommand{\tmop}[1]{\ensuremath{\operatorname{#1}}}

\newtheorem{theorem}{Theorem}[section]

\numberwithin{equation}{section}

 \makeatletter
 
 \newcommand{\Rmnum}[1]{\expandafter\@slowromancap\romannumeral #1@}
 \makeatother
\begin{document}

\title{On the number of orthogonal systems in vector spaces over finite
fields}\author{Le Anh Vinh\\
Mathematics Department\\
Harvard University\\
Cambridge, MA 02138, US\\
vinh@math.harvard.edu}\maketitle

\begin{abstract}
Iosevich and Senger (2008) showed that if a subset of the $d$-dimensional vector space over a finite field is large enough, then it contains many $k$-tuples of mutually orthogonal vectors. In this note, we provide a graph theoretic proof of this result.
\end{abstract}

\section{Introduction}

A classical set of problems in combinatorial geometry deals with the question
of whether a sufficiently large subset of $\mathbbm{R}^d$, $\mathbbm{Z}^d$ or
$\mathbbm{F}_q^d$ contains a given geometric configuration. In a recent paper
\cite{iosevich-senger}, Iosevich and Senger showed that a sufficiently large
subset of $\mathbbm{F}_q^d$, the $d$-dimensional vector space over the finite
field with $q$ elements, contains many $k$-tuple of mutually orthogonal
vectors. Using geometric and character sum machinery, they proved the
following result (see \cite{iosevich-senger} for the motivation of this
result).

\begin{theorem} \label{old}
  (\cite{iosevich-senger}) Let $E \subset \mathbbm{F}_q^d$, such that
  \begin{equation}
    |E| \geqslant C q^{d \frac{k - 1}{k} + \frac{k - 1}{2} + \frac{1}{k}}
  \end{equation}
  with a sufficiently large constant $C > 0$, where $0 < \binom{k}{2} < d$.
  Let $\lambda_k$ be the number of $k$-tuples of $k$ mutually orthogonal
  vectors in $E$. Then
  \begin{equation}
    \lambda_k = (1 + o (1)) \frac{|E|^k}{k!} q^{- \binom{k}{2}} .
  \end{equation}
\end{theorem}

In this note, we provide a different proof to this result using graph
theoretic methods. The main result of this note is the following.

\begin{theorem}\label{main}
  Let $E \subset \mathbbm{F}_q^d$, such that
  \begin{equation}
    |E| \gg q^{\frac{d}{2} + k-1},
  \end{equation}
  where $d \geqslant 2k-1$. Then the number of $k$-tuples of $k$ mutually
  orthogonal vectors in $E$ is
  \begin{equation}
    (1 + o (1)) \frac{|E|^k}{k!} q^{- \binom{k}{2}} .
  \end{equation}
\end{theorem}

Note that Theorem \ref{old} only works in the range $d \geqslant \binom{k}{2}$ (as
larger tuples of mutually orthogonal vectors are out of range of the methods
uses) while Theorem \ref{main} works in a wider range $d \geqslant 2k-1$. Moreover,
Theorem \ref{main} is stronger than Theorem \ref{old} in the same range.

\section{Proof of Theorem \ref{main}}

We call a graph $G = (V, E)$ $(n, d, \lambda)$-graph if $G$ is a $d$-regular
graph on $n$ vertices with the absolute values of each of its eigenvalues but the
largest one is at most $\lambda$. It is well-known that if $\lambda \ll d$
then an $(n, d, \lambda)$-graph behaves similarly as a random graph $G_{n, d /
n}$.  Let $H$ be a fixed graph of order $s$ with $r$ edges and with automorphism
group $\tmop{Aut} (H)$. Using the second moment method, it is not difficult to show that for every constant $p$ the
random graph $G (n, p)$ contains
\begin{equation}
  (1 + o (1)) p^r (1 - p)^{(^s_2) - r} \frac{n^s}{| \tmop{Aut} (H) |}
\end{equation}
induced copies of $H$. Alon extended this result to $(n, d, \lambda)$-graphs.
He proved that every large subset of the set of vertices of an $(n, d,
\lambda)$-graph contains the ``correct'' number of copies of any fixed small
subgraph (Theorem 4.10 in \cite{krivelevich-sudakov}).

\begin{theorem}\label{tool}
  (\cite{krivelevich-sudakov}) Let $H$ be a fixed graph with $r$ edges, $s$
  vertices and maximum degree $\Delta$, and let $G = (V, E)$ be an $(n, d,
  \lambda)$-graph, where, say, $d \leqslant 0.9 n$. Let $m < n$ satisfies $m
  \gg \lambda \left( \frac{n}{d} \right)^{\Delta}$. Then, for every subset $U
  \subset V$ of cardinality $m$, the number of (not necessarily induced) copies
  of $H$ in $U$ is
  \begin{equation}
    (1 + o (1)) \frac{m^s}{| \tmop{Aut} (H) |} \left( \frac{d}{n} \right)^r .
  \end{equation}
\end{theorem}

Note that the above theorem is stated for simple graphs in \cite{krivelevich-sudakov} but there is no different in the proof if we allow loops in the graph $G$. 

We recall a well-known construction of Alon and Krivelevich
\cite{alon}. Let $P G (q, d)$ denote the projective geometry of dimension $d
- 1$ over finite field $\mathbbm{F}_q$. The vertices of $P G (q, d)$
correspond to the equivalence classes of the set of all non-zero vectors $x =
(x_1, \ldots, x_d)$ over $\mathbbm{F}_q$, where two vectors are equivalent if
one is a multiple of the other by an element of the field. Let $G_P (q, d)$
denote the graph whose vertices are the points of $P G (q, d)$ and two (not
necessarily distinct) vertices $x$ and $y$ are adjacent if and only if $x_1 y_1
+ \ldots + x_d y_d = 0$. This construction is well known. In the case $d = 2$,
this graph is called the Erd\"os-R\'enyi graph. It is easy to see that the number
of vertices of $G_P(q,d)$ is $n_{q, d} = (q^d - 1) / (q - 1)$ and that it is $d_{q,
d}$-regular for $d_{q, d} = (q^{d - 1} - 1) / (q - 1)$. The eigenvalues of $G$
are easy to compute (\cite{alon}). Let $A$ be the adjacency matrix of $G$.
Then, by properties of $P G (q, d)$, $A^2 = A A^T = \mu J + (d_{q, d} - \mu)
I$, where $\mu = (q^{d - 2} - 1) / (q - 1)$, $J$ is the all one matrix and $I$
is the identity matrix, both of size $n_{q, d} \times n_{q, d} .$ Thus the
largest eigenvalue of $A$ is $d_{q, d}$ and the absolute value of all other
eigenvalues is $\sqrt{d_{q, d} - \mu} = q^{(d - 2) / 2}$.

Now we are ready to give a proof of Theorem \ref{main}. Let $G (q, d)$ denote the graph whose vertices are the points of
$\mathbbm{F}_q^d - (0, \ldots, 0)$ and two (not necessarily distinct) vertices
$x$ and $y$ are adjacent if and only if they are orthogonal, i.e. $x_1 y_1 +
\ldots + x_d y_d = 0$. Then $G (q, d)$ is just the product of $q - 1$ copies
of $G_P (q, d)$. Therefore, it is easy to see that the number of vertices of
$G$ is $N_{q, d} = (q-1)n_{q,d} = q^d - 1$ and that it is $D_{q, d}$-regular for $D_{q, d} = (q-1)d_{q,d} =
q^{d - 1} - 1$. The eigenvalues of $G (q, d)$ are also easy to compute. Let
$V$ be the adjacency matrix of $G (q, d)$. Then by the properties of $P G (q,
d)$,
\begin{equation}
  V^2 = V V^T = \rho J_{N_{q, d}} + (D_{q, d} - \rho) \bigoplus_{n_{q, d}}
  J_{q - 1},
\end{equation}
where $\rho = (q-1)\mu = q^{d - 2} - 1$, $J_{N_{q, d}}$ is the all one matrix of size
$N_{q, d} \times N_{q, d}$ and $J_{q - 1}$ is the all one matrix of size $(q -
1) \times (q - 1)$. Thus, all eigenvalues of $V^2$ are all eigenvalues of
$(q-1)\rho J_{n_{q, d}} + (q-1)(D_{q, d} - \rho) I_{n_{q,d}}$ and zeros (with $J_{n_{q, d}}$ is
the all one matrix and $I_{n_{q,d}}$ is the identity matrix, both of size $n_{q, d}
\times n_{q, d}$). Therefore, the largest eigenvalue of $V$ is $D_{q, d}$ and
the absolute values of all other eigenvalues are either $\sqrt{(q-1)(D_{q, d} -
\rho)} = (q-1)q^{(d-2)/2}$ or $0$. This implies that $G (q, d)$ is a $(q^d - 1, q^{d - 1} - 1, (q-1)q^{(d-2)/2}$)-graph. Theorem \ref{main} now follows immediately from Theorem \ref{tool}.

\section*{Acknowledgments}

The research is performed during the author's visit at the Erwin Schr\"odinger International Institute for Mathematical Physics. The author would like to thank the ESI for hospitality and financial support during his visit.

\end{document}